\documentclass[submission]{FPSAC2026arxiv}

\newtheorem{thm}{Theorem}

\newtheorem{prop}[thm]{Proposition}
\newtheorem{cor}[thm]{Corollary}

\newtheorem{innercustomthm}{Theorem}
\newenvironment{customthm}[1]
  {\renewcommand\theinnercustomthm{#1}\innercustomthm}
  {\endinnercustomthm}

\theoremstyle{definition}
\newtheorem{defn}[thm]{Definition}
\newtheorem{example}[thm]{Example}
\newtheorem{remark}[thm]{Remark}

\numberwithin{thm}{section}

\usepackage{tikz}
\usetikzlibrary{calc, math}

\DeclareMathOperator{\Gr}{Gr}
\DeclareMathOperator{\init}{in}
\DeclareMathOperator{\Spec}{Spec}

\DeclareMathOperator{\Trop}{Trop}

\DeclareMathOperator{\sgn}{sgn}

\newcommand{\la}{\lambda}

\newcommand{\cM}{\mathcal M}

\newcommand{\cA}{\mathcal A}

\newcommand{\cT}{\mathcal T}
\newcommand{\bC}{\mathbb{C}}

\newcommand{\bR}{\mathbb{R}}

\newcommand{\rA}{\mathrm A}
\newcommand{\rC}{\mathrm C}

\newcommand{\bs}{\backslash}

\newcommand{\mr}{\mathring}
\newcommand{\ol}{\overline}
\newcommand{\ph}{\varphi}

\newcommand{\vv}{\mathbf v}

\mathcode`l="8000
\begingroup
\makeatletter
\lccode`\~=`\l
\DeclareMathSymbol{\lsb@l}{\mathalpha}{letters}{`l}
\lowercase{\gdef~{\ifnum\the\mathgroup=\m@ne \ell \else \lsb@l \fi}}%
\endgroup

\title[Tropical cluster varieties and phylogenetic trees]{Tropical cluster varieties, phylogenetic trees, and generalized associahedra}

\author{Igor Makhlin\thanks{\href{mailto:iymakhlin@gmail.com}{iymakhlin@gmail.com}. Igor Makhlin was partially supported by SFB-TRR 195 ``Symbolic Tools in Mathematics and their Application''.}\addressmark{1}}

\address{\addressmark{1}  Technische Universität Berlin, Fakultät II Mathematik und Naturwissenschaften, Institut für Mathematik, FG Diskrete Mathematik/Geometrie}

\received{\today}

\abstract{We explicitly describe the tropicalization of a type C cluster variety by identifying it with the space of axially symmetric phylogenetic trees. We also study the signed tropicalizations of this cluster variety, realizing them as subfans of the tropicalization that are dual to either associahedra or cyclohedra.}


\usepackage[backend=bibtex]{biblatex}
\begin{document}

\maketitle

\section*{Introduction}

At least two notions of tropicalization have been widely studied for cluster algebras and cluster varieties. These are \textit{positive tropicalizations}, defined by Speyer and Williams in~\cite{SpeWil2005}, and the related \textit{Fock--Goncharov tropicalizations}, due to Fock and Goncharov~\cite{FocGon2006,FocGon2009} and Gross, Hacking, Keel, and Kontsevich~\cite{GHKK2018}.

However, much less is known about \textit{tropicalizations} of cluster varieties as commonly understood in tropical algebraic geometry today, following~\cite{BieGro1984,Stu2002,Mikh2005,EinKapLin2006} and others.
In broad terms, the tropicalization of an algebraic variety is its ``combinatorial shadow'': a polyhedral complex that captures various geometric properties of the variety, from dimension and degree to symmetry groups and intersection theory.  
In the pioneering work~\cite{SpeStu2004}, Speyer and Sturmfels gave an explicit description of the tropicalization of the Grassmannian $\Gr(2,n)$, realizing it as the space of phylogenetic trees, introduced in~\cite{RobWhi1996,BilHolVog2001}. This offers a direct connection to cluster algebra theory, since the affine cone over $\Gr(2,n)$ arises as the prototypical example of a cluster variety: that of finite type $\rA$.

While highly consequential, the result of~\cite{SpeStu2004} has resisted attempts of extending it to higher Grassmannians.
To the author's knowledge, it also remains the only explicit description of an infinite family of tropicalized cluster varieties. Thus, the broad motivation behind our work is twofold: to obtain further families of tropicalized cluster varieties, and to show that this may be a more natural (or at least more accessible) direction for generalizing the construction of Speyer and Sturmfels. Specifically, we make a first step towards this goal by describing tropicalized cluster varieties of finite type C.

We build on the work~\cite{CoxMakh2025}, where Cox and the author introduced a polyhedral fan termed the \textit{space of axially symmetric phylogenetic trees} (or \textit{ASPTs}). This fan parametrizes phylogenetic trees with leaves labeled by $\{1,\dots,n,\overline 1,\dots,\overline n\}$ that satisfy a natural symmetry condition, cf.\ Figure~\ref{fig:graphn=3}. Our main result is as follows.

\begin{customthm}{A}[Thm.~\ref{thm:main}, Cor.~\ref{cor:fullrank}]
For a full-rank cluster variety $X(\tilde B)$ of finite type $\mathrm C_{n-1}$, the \textit{tropical cluster variety} $\Trop X(\tilde B)$ is, modulo lineality space, equal to the space of ASPTs.   
\end{customthm}

We also present a second result that exhibits a finer structure on the type C cluster variety, given in terms of generalized permutahedra. 
Our result characterizes the \textit{signed tropicalizations} of $X(\tilde B)$. These are certain subfans of $\Trop X(\tilde B)$ that arise as natural generalizations of the positive tropicalization. A signed tropicalization is associated with a \textit{sign pattern}: a sequence of signs, one for every coordinate of the ambient affine space.
\begin{customthm}{B}[Cor.~\ref{cor:signedtrops}]\label{thm:B}
For every sign pattern that occurs in $X(\tilde B)$, the respective signed tropicalization of $X(\tilde B)$ is, modulo lineality space, combinatorially equivalent to either the dual fan of an associahedron or the dual fan of a cyclohedron.
\end{customthm}

Theorem~\ref{thm:B} follows from an explicit description of the signed tropicalizations of the \textit{cluster configuration space} $\cM$, a certain geometric quotient of $X(\tilde B)$. The signed tropicalizations of $\cM$ are enumerated by dihedral orderings that are centrally or axially symmetric, see Theorem~\ref{thm:signedtrops}.


This is an extended abstract of the preprint~\cite{Mak2025}.

\section{ASPTs and the space of ASPTs}

We fix an integer $n\ge 3$ and consider the $2n$-element set $N=\{1,\dots,n,\ol 1,\dots,\ol n\}$. For $i\in[1,n]$ we set $\ol{\ol i}=i$ and $|\ol i|=|i|=i$.

A \textit{phylogenetic tree} is a pair $(T,\vv)$ consisting of a tree $T$ with $2n$ leaves and no vertices of degree 2, together with a bijection $\vv$ from $N$ to the set of leaf vertices in $T$. We identify phylogenetic trees that differ only by label-preserving graph isomorphism.


Let $P$ denote a regular $2n$-gon. We say that two diagonals of $P$ cross if they share exactly one interior point. We use the term \textit{subdivision of $P$} to refer to a set of pairwise non-crossing diagonals. A \textit{labeling} of $P$ is a bijection from $N$ to the set of sides of $P$. Given a subdivision $\Theta$ of $P$ and a labeling $\ph$ of $P$, we obtain a phylogenetic tree $\cT_{\Theta,\ph}=(T,\vv)$ as the dual graph. Explicitly:
\begin{itemize}
\item the non-leaf vertices of $T$ correspond to the polygonal cells formed by the diagonals in $\Theta$, with two such vertices adjacent if and only if the respective cells share a side;
\item for $a\in N$, the leaf $\vv(a)$ is adjacent to the non-leaf vertex that corresponds to the cell containing the side $\ph(a)$.
\end{itemize}
It is clear that any phylogenetic tree has the form $\cT_{\Theta,\ph}$ for an appropriate choice of $\Theta$ and $\ph$. However, for a given phylogenetic tree this choice is not unique, even modulo symmetries of $P$.

\begin{example}
In Figure~\ref{fig:subdivs}, we see three possible subdivisions $\Theta$ for the case $n=4$, the respective dual tree $T$ is shown in red. In each case we have also chosen a labeling $\ph$, and labeled each side $s$ with $\ph^{-1}(s)$. As a result we obtain the phylogenetic tree $\cT_{\Theta,\ph}=(T,\vv)$, with each leaf vertex $v$ labeled by $\vv^{-1}(v)$. Note that in the first and third cases we obtain the same phylogenetic tree.

\begin{figure}[h!tbp]
\centering

\begin{tikzpicture}[scale=2]
\foreach \i in {1,...,8}{
\coordinate (P\i) at ({cos(90-45*(\i-1))},{sin(90-45*(\i-1))});
}

\draw (P1)--(P2)--(P3)--(P4)--(P5)--(P6)--(P7)--(P8)--cycle;
\draw (P1)--(P3) (P1)--(P5) (P1)--(P7);

\coordinate (C1) at (0.571, 0.589);   
\coordinate (C2) at (0.354, 0.146);   
\coordinate (C3) at (-0.354, 0.146);  
\coordinate (C4) at (-0.571, 0.589);  

\foreach \c in {C1,C2,C3,C4}
\fill[red] (\c) circle (1.4pt);

\draw[red] (C1)--(C2) (C3)--(C4) (C2)--(C3);

\foreach \i/\j/\name in {1/2/M12,2/3/M23,3/4/M34,4/5/M45,5/6/M56,6/7/M67,7/8/M78,8/1/M81}{
  \coordinate (\name) at ($ (P\i)!0.5!(P\j) $);
  \coordinate (shifted\name) at ($ (\name) + 0.1*($( \name )-(0,0)$) $);
  \fill[red] (shifted\name) circle (1.4pt);
}

\draw[red] (shiftedM12)--(C1) (shiftedM23)--(C1);
\draw[red] (shiftedM34)--(C2) (shiftedM45)--(C2);
\draw[red] (shiftedM56)--(C3) (shiftedM67)--(C3);
\draw[red] (shiftedM78)--(C4) (shiftedM81)--(C4); 

\node at (shiftedM12) [above right, black] {$2$};
\node at (shiftedM23) [above right, black] {$\ol 4$};
\node at (shiftedM34) [above right, black] {$1$};
\node at (shiftedM45) [below right, black] {$\ol 3$};
\node at (shiftedM56) [below left, black] {$3$};
\node at (shiftedM67) [above left, black] {$\ol 1$};
\node at (shiftedM78) [above left, black] {$4$};
\node at (shiftedM81) [above left, black] {$\ol 2$};

\path (P1) -- (P5) node[pos=0.25, left] {$\delta_0$};
\end{tikzpicture}
\hspace{5mm}
\begin{tikzpicture}[scale=2]
\foreach \i in {1,...,8}{
  \coordinate (P\i) at ({cos(90-45*(\i-1))},{sin(90-45*(\i-1))});
}

\draw (P1)--(P2)--(P3)--(P4)--(P5)--(P6)--(P7)--(P8)--cycle;
\draw (P1)--(P4) (P1)--(P6) (P4)--(P6);

\coordinate (C1) at (0.55,0.25);   
\coordinate (C2) at (0,0.2);      
\coordinate (C3) at (-0.55,0.25);  
\coordinate (C4) at (0,-0.85);     

\foreach \c in {C1,C2,C3,C4}
  \fill[red] (\c) circle (1.4pt);

\draw[red] (C2)--(C1);
\draw[red] (C2)--(C3);
\draw[red] (C2)--(C4);

\foreach \i/\j/\name in {1/2/M12,2/3/M23,3/4/M34,4/5/M45,5/6/M56,6/7/M67,7/8/M78,8/1/M81}{
  \coordinate (\name) at ($ (P\i)!0.5!(P\j) $);
  \coordinate (shifted\name) at ($ (\name) + 0.1*($( \name )-(0,0)$) $);
  \fill[red] (shifted\name) circle (1.4pt);
}

\draw[red] (shiftedM12)--(C1);
\draw[red] (shiftedM23)--(C1);
\draw[red] (shiftedM34)--(C1);
\draw[red] (shiftedM45)--(C4);
\draw[red] (shiftedM56)--(C4);
\draw[red] (shiftedM67)--(C3);
\draw[red] (shiftedM78)--(C3);
\draw[red] (shiftedM81)--(C3);

\node at (shiftedM12) [above right, black] {$1$};
\node at (shiftedM23) [above right, black] {$\ol 3$};
\node at (shiftedM34) [above right, black] {$2$};
\node at (shiftedM45) [below right, black] {$4$};
\node at (shiftedM56) [below left, black] {$\ol 4$};
\node at (shiftedM67) [above left, black] {$\ol 2$};
\node at (shiftedM78) [above left, black] {$3$};
\node at (shiftedM81) [above left, black] {$\ol 1$};

\draw[dotted] (P1) -- (P5) node[pos=0.28, left, xshift=0.1cm] {$\delta_0$};
\end{tikzpicture}
\hspace{5mm}
\begin{tikzpicture}[scale=2]
\foreach \i in {1,...,8}{
  \coordinate (P\i) at ({cos(90-45*(\i-1))},{sin(90-45*(\i-1))});
}

\draw (P1)--(P2)--(P3)--(P4)--(P5)--(P6)--(P7)--(P8)--cycle;
\draw (P2)--(P6) (P2)--(P8) (P4)--(P6);

\coordinate (C1) at (0,0.8);     
\coordinate (C2) at (0.38,-0.25);  
\coordinate (C3) at (0,-0.8);    
\coordinate (C4) at (-0.38,0.25); 

\foreach \c in {C1,C2,C3,C4}
  \fill[red] (\c) circle (1.4pt);

\draw[red] (C1)--(C4)--(C2)--(C3);

\foreach \i/\j/\name in {1/2/M12,2/3/M23,3/4/M34,4/5/M45,5/6/M56,6/7/M67,7/8/M78,8/1/M81}{
  \coordinate (\name) at ($ (P\i)!0.5!(P\j) $);
  \coordinate (shifted\name) at ($ (\name) + 0.1*($( \name )-(0,0)$) $);
  \fill[red] (shifted\name) circle (1.4pt);
}

\draw[red] (shiftedM12)--(C1);
\draw[red] (shiftedM23)--(C2);
\draw[red] (shiftedM34)--(C2);
\draw[red] (shiftedM45)--(C3);
\draw[red] (shiftedM56)--(C3);
\draw[red] (shiftedM67)--(C4);
\draw[red] (shiftedM78)--(C4);
\draw[red] (shiftedM81)--(C1);

\node at (shiftedM12) [above right, black] {$4$};
\node at (shiftedM23) [above right, black] {$1$};
\node at (shiftedM34) [above right, black] {$\ol 3$};
\node at (shiftedM45) [below right, black] {$2$};
\node at (shiftedM56) [below left, black] {$\ol 4$};
\node at (shiftedM67) [above left, black] {$\ol 1$};
\node at (shiftedM78) [above left, black] {$3$};
\node at (shiftedM81) [above left, black] {$\ol 2$};
\end{tikzpicture}

\vspace{-1cm}
\caption{}
\label{fig:subdivs}
\end{figure}    
\end{example}

Now, we fix one of the $n$ longest diagonals of $P$ and denote it by $\delta_0$. We say that a subdivision of $P$ is \textit{axially symmetric} if with every diagonal it also contains its mirror image with respect to $\delta_0$. Such a subdivision can only contain diagonals that do not cross $\delta_0$ or are perpendicular to $\delta_0$. A labeling of $P$ is \textit{axially symmetric} if the sides $\ph(i)$ and $\ph(\ol i)$ are symmetric to each other with respect to $\delta_0$ for every $i\in [1,n]$. 

\begin{defn}
An \textit{axially symmetric phylogenetic tree} (or \textit{ASPT}) is a phylogenetic tree that has the form $\cT_{\Theta,\ph}$ for an axially symmetric subdivision $\Theta$ and an axially symmetric labeling $\ph$.
\end{defn}

\begin{example}
In the first and second pictures in Figure~\ref{fig:subdivs}, both $\Theta$ and $\ph$ are axially symmetric, and we obtain ASPTs. In the third picture, neither $\Theta$ nor $\ph$ is axially symmetric but $\cT_{\Theta,\ph}$ is still an ASPT. A non-ASPT can be obtained from any of the three by switching labels 1 and 4.    
\end{example} 

\begin{example}\label{ex:ASPTn=3}
For $n=3$, every ASPT has one of the 7 forms shown in Figure~\ref{fig:ASPTn=3}, where the set $\{a,b,c\}\subset N$ satisfies $\{|a|,|b|,|c|\}=\{1,2,3\}$.

\begin{figure}[h!tbp]

\centering

\hspace{3mm}
\begin{tikzpicture}[x=5mm, y=5mm]
\node at (-2,2) {1.};

\node[circle, fill=black, inner sep=2] (p) at (0,0) {};

\node[circle, fill=black, inner sep=2] (a) at (-1,1) {};
\node[circle, fill=black, inner sep=2] (b) at (1,1) {};
\node[circle, fill=black, inner sep=2] (c) at (-2,0) {};
\node[circle, fill=black, inner sep=2] (d) at (2,0) {};
\node[circle, fill=black, inner sep=2] (e) at (-1,-1) {};
\node[circle, fill=black, inner sep=2] (f) at (1,-1) {};

\node at (-1.5,1) {$a$};
\node at (1.5,1) {$\ol a$};
\node at (-2.5,0) {$b$};
\node at (2.5,0) {$\ol b$};
\node at (-1.5,-1) {$c$};
\node at (1.5,-1) {$\ol c$};

\draw (p) -- (a)
      (p) -- (b)
      (p) -- (c)
      (p) -- (d)
      (p) -- (e)
      (p) -- (f);
\end{tikzpicture} 
\hspace{3mm}
\begin{tikzpicture}[x=5mm, y=5mm]
\node at (-2.5,2) {2.};

\node[circle, fill=black, inner sep=2] (p) at (0,1) {};
\node[circle, fill=black, inner sep=2] (q) at (0,0) {};

\node[circle, fill=black, inner sep=2] (a) at (-1,2) {};
\node[circle, fill=black, inner sep=2] (b) at (1,2) {};
\node[circle, fill=black, inner sep=2] (c) at (-2,0) {};
\node[circle, fill=black, inner sep=2] (d) at (2,0) {};
\node[circle, fill=black, inner sep=2] (e) at (-1,-1) {};
\node[circle, fill=black, inner sep=2] (f) at (1,-1) {};

\node at (-1.5,2) {$a$};
\node at (1.5,2) {$\ol a$};
\node at (-2.5,0) {$b$};
\node at (2.5,0) {$\ol b$};
\node at (-1.5,-1) {$c$};
\node at (1.5,-1) {$\ol c$};

\draw (p) -- (a)
      (p) -- (b)
      (q) -- (c)
      (q) -- (d)
      (q) -- (e)
      (q) -- (f)
      (p) -- (q);
\end{tikzpicture} 
\hspace{3mm}
\begin{tikzpicture}[x=5mm, y=5mm]
\node at (-2,2) {3.};

\node[circle, fill=black, inner sep=2] (p) at (0,0) {};
\node[circle, fill=black, inner sep=2] (q) at (1,0) {};

\node[circle, fill=black, inner sep=2] (a) at (-1,1) {};
\node[circle, fill=black, inner sep=2] (b) at (2,1) {};
\node[circle, fill=black, inner sep=2] (c) at (-2,0) {};
\node[circle, fill=black, inner sep=2] (d) at (3,0) {};
\node[circle, fill=black, inner sep=2] (e) at (-1,-1) {};
\node[circle, fill=black, inner sep=2] (f) at (2,-1) {};

\node at (-1.5,1) {$a$};
\node at (2.5,1) {$\ol a$};
\node at (-2.5,0) {$b$};
\node at (3.5,0) {$\ol b$};
\node at (-1.5,-1) {$c$};
\node at (2.5,-1) {$\ol c$};

\draw (p) -- (a)
      (q) -- (b)
      (p) -- (c)
      (q) -- (d)
      (p) -- (e)
      (q) -- (f)
      (p) -- (q);
\end{tikzpicture} 
\hspace{3mm}
\begin{tikzpicture}[x=5mm, y=5mm]
\node at (-2.5,2) {4.};

\node[circle, fill=black, inner sep=2] (p) at (0,1) {};
\node[circle, fill=black, inner sep=2] (q) at (0,0) {};
\node[circle, fill=black, inner sep=2] (r) at (0,-1) {};

\node[circle, fill=black, inner sep=2] (a) at (-1,2) {};
\node[circle, fill=black, inner sep=2] (b) at (1,2) {};
\node[circle, fill=black, inner sep=2] (c) at (-2,0) {};
\node[circle, fill=black, inner sep=2] (d) at (2,0) {};
\node[circle, fill=black, inner sep=2] (e) at (-1,-2) {};
\node[circle, fill=black, inner sep=2] (f) at (1,-2) {};

\node at (-1.5,2) {$a$};
\node at (1.5,2) {$\ol a$};
\node at (-2.5,0) {$b$};
\node at (2.5,0) {$\ol b$};
\node at (-1.5,-2) {$c$};
\node at (1.5,-2) {$\ol c$};

\draw (p) -- (a)
      (p) -- (b)
      (q) -- (c)
      (q) -- (d)
      (r) -- (e)
      (r) -- (f)
      (p) -- (q)
      (q) -- (r);
\end{tikzpicture}
\vspace{2mm}

\hspace{3mm}
\begin{tikzpicture}[x=5mm, y=5mm]
\node at (-2.5,1.5) {5.};

\node[circle, fill=black, inner sep=2] (p) at (0,0) {};
\node[circle, fill=black, inner sep=2] (q) at (-1,-1) {};
\node[circle, fill=black, inner sep=2] (r) at (1,-1) {};

\node[circle, fill=black, inner sep=2] (a) at (-1,1) {};
\node[circle, fill=black, inner sep=2] (b) at (1,1) {};
\node[circle, fill=black, inner sep=2] (c) at (-2,-1) {};
\node[circle, fill=black, inner sep=2] (d) at (2,-1) {};
\node[circle, fill=black, inner sep=2] (e) at (-1,-2) {};
\node[circle, fill=black, inner sep=2] (f) at (1,-2) {};

\node at (-1.5,1) {$a$};
\node at (1.5,1) {$\ol a$};
\node at (-2.5,-1) {$b$};
\node at (2.5,-1) {$\ol b$};
\node at (-1.5,-2) {$c$};
\node at (1.5,-2) {$\ol c$};

\draw (p) -- (a)
      (p) -- (b)
      (q) -- (c)
      (r) -- (d)
      (q) -- (e)
      (r) -- (f)
      (p) -- (q)
      (p) -- (r);
\end{tikzpicture}
\hspace{3mm}
\begin{tikzpicture}[x=5mm, y=5mm]
\node at (-2.5,2) {6.};

\node[circle, fill=black, inner sep=2] (p) at (0,1) {};
\node[circle, fill=black, inner sep=2] (q) at (0,0) {};
\node[circle, fill=black, inner sep=2] (r) at (-1,-1) {};
\node[circle, fill=black, inner sep=2] (s) at (1,-1) {};

\node[circle, fill=black, inner sep=2] (a) at (-1,2) {};
\node[circle, fill=black, inner sep=2] (b) at (1,2) {};
\node[circle, fill=black, inner sep=2] (c) at (-2,-1) {};
\node[circle, fill=black, inner sep=2] (d) at (2,-1) {};
\node[circle, fill=black, inner sep=2] (e) at (-1,-2) {};
\node[circle, fill=black, inner sep=2] (f) at (1,-2) {};

\node at (-1.5,2) {$a$};
\node at (1.5,2) {$\ol a$};
\node at (-2.5,-1) {$b$};
\node at (2.5,-1) {$\ol b$};
\node at (-1.5,-2) {$c$};
\node at (1.5,-2) {$\ol c$};

\draw (p) -- (a)
      (p) -- (b)
      (r) -- (c)
      (s) -- (d)
      (r) -- (e)
      (s) -- (f)
      (p) -- (q)
      (q) -- (r)
      (q) -- (s);
\end{tikzpicture} 
\hspace{3mm}
\begin{tikzpicture}[x=5mm, y=5mm]
\node at (-2.5,1.5) {7.};

\node[circle, fill=black, inner sep=2] (p) at (0,0) {};
\node[circle, fill=black, inner sep=2] (q) at (1,0) {};
\node[circle, fill=black, inner sep=2] (r) at (-1,-1) {};
\node[circle, fill=black, inner sep=2] (s) at (2,-1) {};

\node[circle, fill=black, inner sep=2] (a) at (-1,1) {};
\node[circle, fill=black, inner sep=2] (b) at (2,1) {};
\node[circle, fill=black, inner sep=2] (c) at (-2,-1) {};
\node[circle, fill=black, inner sep=2] (d) at (3,-1) {};
\node[circle, fill=black, inner sep=2] (e) at (-1,-2) {};
\node[circle, fill=black, inner sep=2] (f) at (2,-2) {};

\node at (-1.5,1) {$a$};
\node at (2.5,1) {$\ol a$};
\node at (-2.5,-1) {$b$};
\node at (3.5,-1) {$\ol b$};
\node at (-1.5,-2) {$c$};
\node at (2.5,-2) {$\ol c$};

\draw (p) -- (a)
      (q) -- (b)
      (r) -- (c)
      (s) -- (d)
      (r) -- (e)
      (s) -- (f)
      (p) -- (q)
      (p) -- (r)
      (q) -- (s);
\end{tikzpicture}

\vspace{-1cm}
\caption{}\label{fig:ASPTn=3}

\end{figure}
\end{example}

It is clear from the definition that every ASPT $(T,\vv)$ is equipped with a unique involution $\sigma$ of $T$ that exchanges the vertices $\vv(i)$ and $\vv(\ol i)$. We will refer to $\sigma$ as the \textit{symmetry} of $(T,\vv)$.

\begin{defn}
A \textit{weighted phylogenetic tree} $(T,\vv,l)$ is a phylogenetic tree $(T,\vv)$ together with a weight function $l$ from the edge set of $T$ to $\bR$ such that $l(e)>0$ for every non-leaf edge $e$ (but not necessarily for the leaf edges). We also say that $l$ is a \textit{weighting} of $(T,\vv)$.
\end{defn}

A weighted phylogenetic tree $(T,\vv,l)$ defines a ``distance'' function $d_{T,\vv,l}\colon N^2\to\bR$. Specifically, $d_{T,\vv,l}(a,b)$ is equal to the total of $l(e)$ over all edges $e$ lying in the non-self-intersecting path between $\vv(a)$ and $\vv(b)$.
The function $d_{T,\vv,l}(a,b)$ is symmetric and satisfies $d_{T,\vv,l}(a,a)=0$. 

\begin{defn}
An \textit{axially symmetric weighted phylogenetic tree} (or \textit{ASWPT}) is a weighted phylogenetic tree $(T,\vv,l)$ such that $(T,\vv)$ is an ASPT and $l(\sigma(e))=l(e)$ for every edge $e$, where $\sigma$ is the symmetry of $(T,\vv)$.
\end{defn}

Evidently, an ASWPT $(T,\vv,l)$ has the property that $d_{T,\vv,l}(a,b)=d_{T,\vv,l}(\ol a,\ol b)$. This leads us to consider the set of all pairs $\{a,b\}\subset N$ modulo the equivalence relation $\{a,b\}\sim\{\ol a,\ol b\}$. Such equivalence classes are enumerated by the following $n^2$-element set $D\subset N^2$: 
\[D=\{(i,j)\mid 1\le i<j\le n\}\sqcup\{(i,\ol j)\mid 1\le i\le j\le n\},\]
since every class contains exactly one pair $\{a,b\}$ with $(a,b)\in D$.

Now we can associate a point $w(T,\vv,l)\in\bR^D$ with every ASWPT $(T,\vv,l)$ by setting 
\[w(T,\vv,l)_{a,b}=d_{T,\vv,l}(a,b).\] 
Subsequently, for an ASPT $(T,\vv)$ we consider the set
\[C_{T,\vv}=\{w(T,\vv,l)\,|\,l\text{ is a weighting of }(T,\vv)\}\subset\bR^D.\]



Before listing the key properties of the sets $C_{T,\vv}$, we introduce one more notation. For an ASPT $(T,\vv)$ with symmetry $\sigma$, we let $k(T,\vv)$ denote the number of $\sigma$-orbits in the edge set of $T$ (with each orbit containing 1 or 2 elements).

\begin{prop}\label{prop:spaceofASPTs}
\hfill
\begin{enumerate}[label=(\alph*)]
\item If ASWPTs $(T_1,\vv_1,l_1)$ and $(T_2,\vv_2,l_2)$ are distinct, then $w(T_1,\vv_1,l_1)\neq w(T_2,\vv_2,l_2)$. If ASPTs $(T_1,\vv_1)$ and $(T_2,\vv_2)$ are distinct, then $C_{T_1,\vv_1}$ and $C_{T_2,\vv_2}$ are disjoint.
\item For an ASPT $(T,\vv)$, the set $C_{T,\vv}$ is a relatively open polyhedral cone of dimension $k(T,\vv)$: the product of a simplicial cone of dimension $k(T,\vv)-n$ and an $n$-dimensional affine space.
\item The collection of cones $C_{T,\vv}$ with $(T,\vv)$ ranging over all ASPTs is a polyhedral fan in $\bR^D$ of pure dimension $2n-1$ and with a lineality space of dimension $n$.
\end{enumerate}    
\end{prop}

We refer to the polyhedral fan part (c) of the proposition as the \textit{space of ASPTs}, and denote this fan by $\Omega$.

\begin{example}
For $n=3$, the space of ASPTs has dimension $5$ and a lineality space of dimension 3. Thus, combinatorially, it is the product of $\bR^3$ and a fan over a graph. This graph is shown in Figure~\ref{fig:graphn=3}. It has 13 vertices (corresponding to 4-dimensional cones in the fan, enumerated by ASPTs of the forms 2, 3 and 5 in Example~\ref{ex:ASPTn=3}) and 21 edges (corresponding to maximal cones, enumerated by ASPTs of the forms 4, 6 and 7). Some elements of the graph are labeled by the respective ASPTs. 

One may also notice that the ASPT labeling the red vertex is obtained from both of the ASPTs labeling incident red edges by contracting all edges (of the ASPT) that lie in a single $\sigma$-orbit. This is a general rule that specifies when one cone in $\Omega$ arises as a facet of another.

\begin{figure}[h!tbp]
    \centering
    \input{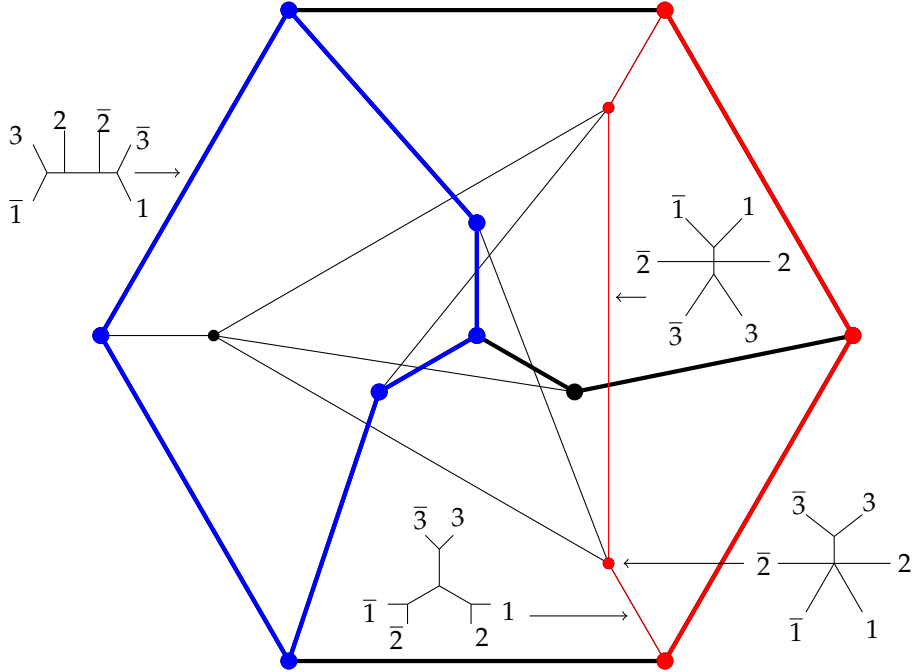}
    \tikzset{hexagon/.pic={
    
}}

\tikzset{
    pics/chexagon/.style args={#1/#2/#3/#4/#5/#6/#7}{
      code = {
        \foreach \i in {0, 1, 2, 3, 4, 5} {
            \coordinate (k\i) at (-60*\i+120:1);
            \coordinate (p\i) at ($.6*(-60*\i+120:1) + .6*(-60*\i+180:1)$);
            
        }

        \node at ($1.15*(p1)$){\scriptsize$#1$};
        \node at ($1.15*(p2)$){\scriptsize$#2$};
        \node at ($1.15*(p3)$){\scriptsize$#3$};
        \node at ($1.15*(p4)$){\scriptsize$-#1$};
        \node at ($1.15*(p5)$){\scriptsize$-#2$};
        \node at ($1.15*(p0)$){\scriptsize$-#3$};

        \draw (k0)--(k1)--(k2)--(k3)--(k4)--(k5)--cycle;

        \draw (k#4)--(k#5);
        \tikzmath{\n = int(mod(#4 + 3, 6)); \m = int(mod(#5 + 3, 6));}
        \draw (k\n)--(k\m);

        \draw (k#6)--(k#7);
        \tikzmath{\n = int(mod(#6 + 3, 6)); \m = int(mod(#7 + 3, 6));}
        \draw (k\n)--(k\m);
    }}}

\tikzset{
    pics/ahexagon/.style args={#1/#2/#3/#4/#5/#6/#7}{
      code = {
        \foreach \i in {0, 1, 2, 3, 4, 5} {
            \coordinate (h\i) at (-60*\i+120:1);
            \coordinate (j\i) at ($.6*(-60*\i+120:1) + .6*(-60*\i+180:1)$);
            
        }

        \node at ($1.15*(j1)$){\scriptsize$#1$};
        \node at ($1.15*(j2)$){\scriptsize$#2$};
        \node at ($1.15*(j3)$){\scriptsize$#3$};
        \node at ($1.15*(j4)$){\scriptsize$-#3$};
        \node at ($1.15*(j5)$){\scriptsize$-#2$};
        \node at ($1.15*(j0)$){\scriptsize$-#1$};

        \draw (h0)--(h1)--(h2)--(h3)--(h4)--(h5)--cycle;

        \draw (h#4)--(h#5);
        \tikzmath{\n = int(mod(#4, 6)); \m = int(mod(-(#5 - 4) + 4, 6));}
        \draw (h\n)--(h\m);

        \draw (h#6)--(h#7);
        \tikzmath{\n = int(mod(#6, 6)); \m = int(mod(-(#5 - 4) + 4, 6));}
        \draw (h\n)--(h\m);
    }}}

{\begin{tikzpicture}[scale=2]
    \foreach \i in {0, 1, 2, 3, 4, 5} {
        \coordinate (a\i) at (60*\i:2.5);
        \coordinate (b\i) at (60*\i:1.75);
        \coordinate (c\i) at (60*\i+30:.75);
        
        \coordinate (l\i) at ($.55*(60*\i:3) + .55*(60*\i+60:3)$);
    }
    \coordinate (alpha) at (0,0);

    \draw[ultra thick] (a0)--(a1)--(a2)--(a3)--(a4)--(a5)--cycle;
    \draw (b1)--(b3)--(b5)--cycle;
    \draw[ultra thick] (a2)--(c1);
    \draw (c1)--(b5);
    \draw[ultra thick] (a4)--(c3);
    \draw (c3)--(b1);
    \draw[ultra thick] (a0)--(c5);
    \draw (c5)--(b3);
    \draw[ultra thick] (c1)--(alpha);
    \draw[ultra thick] (c3)--(alpha);
    \draw[ultra thick] (c5)--(alpha);
    \draw (a1)--(b1) (a3)--(b3) (a5)--(b5);






    \draw[ultra thick, red] (a5)--(a0)--(a1);
    \draw[red] (a1)--(b1)--(b5)--(a5);

    \draw[ultra thick, blue] (alpha)--(c3)--(a4)--(a3)--(a2)--(c1)--cycle;
    
    \foreach \i in {0,1,2,3,4,5} {
        \filldraw (a\i) circle (1.5pt);

        \tikzmath{\n = int(mod(\i, 2));}
        \ifnum \n = 1
        \filldraw (b\i) circle (1pt);
        \filldraw (c\i) circle (1.5pt);
        \fi
    }
    \filldraw (alpha) circle (1.5pt);
    
    \filldraw[red] (b1) circle (1pt);
    \filldraw[red] (b5) circle (1pt);
    \filldraw[red] (a0) circle (1.5pt);
    \filldraw[red] (a1) circle (1.5pt);
    \filldraw[red] (a5) circle (1.5pt);

    \filldraw[blue] (a2) circle (1.5pt);
    \filldraw[blue] (a3) circle (1.5pt);
    \filldraw[blue] (a4) circle (1.5pt);
    \filldraw[blue] (c1) circle (1.5pt);
    \filldraw[blue] (c3) circle (1.5pt);
    \filldraw[blue] (alpha) circle (1.5pt);

    \coordinate (T1) at ($(b1)!.35!(b5)$);
    \pic at ($(T1) + (.7,0)$) [scale=.75] {dtree=\overline{2}/\overline 1/1/2/3/\overline{3}};
    \draw[->] ($(T1) + (.25,-.2)$)--($(T1) + (.05,-.2)$);

    \coordinate (T2) at ($(a2)!.5!(a3)$);
    \pic at ($(T2) + (-.75,0)$) [scale=.75] {etree=\overline{2}/\overline{3}/1/2/\overline{1}/3};
    \draw[->] ($(T2) + (-.4,0)$)--($(T2) + (-.1,0)$);

    \coordinate (T3) at ($(a4)!.4!(a5)$);
    \pic at ($(T3) + (0,.5)$) [scale=.75] {ftree=\overline{3}/3/1/2/\overline{2}/\overline{1}};
    \draw[->] ($(T3) + (.6,.3)$)--($(T3) + (1.25,.3)$);

    \pic at ($(b5) + (1.5,0)$) [scale=.75] {btree=\ol 2/\ol 3/3/2/1/\ol 1};
    \draw[->] ($(b5) + (.9,0)$)--($(b5) + (.1,0)$);
    
\end{tikzpicture}}    
    \caption{For $n=3$, the space of ASPTs is, modulo lineality space, a fan over this graph. Some elements are labeled by the respective ASPTs.}
    \label{fig:graphn=3}
\end{figure}
\end{example}

\section{Type C cluster algebras and the main theorem}

Our approach is to first consider the \textit{special cluster algebra of type C} discussed by Fomin and Zelevinsky in~\cite[\S12.3]{FomZel2003}. This choice of coefficients allows for particularly natural and explicit forms of our results. We will then show how these extend to arbitrary full rank geometric type cluster algebras of the same finite type. 

We denote the chosen cluster algebra of type $\rC_{n-1}$ by $\cA$. Rather than giving the original definition of $\cA$ in terms of mutations, we recall a realization of $\cA$ as of a subalgebra in a polynomial ring, also due to~\cite{FomZel2003}.


As a $\bC$-algebra, $\cA$ is generated by $n^2$ elements $\Delta_{a,b}$, $(a,b)\in D$. These elements form the set of all cluster variables and frozen variables. Let $S$ denote the polynomial ring $\bC[z_{t,i}]_{t\in\{1,2\},\, i\in[1,n]}$ in $2n$ variables. We use the following as a definition of $\cA$.

\begin{prop}[{\cite[Proposition 12.13]{FomZel2003}}]\label{prop:AinS1}
There is an injective homomorphism $\iota$ from $\cA$ to $S$ defined on the generators by
\begin{align*}
\iota(\Delta_{i,j}) & = z_{1,i}z_{2,j}-z_{1,j}z_{2,i},\quad 1\le i<j\le n,\\
\iota(\Delta_{i,\ol j}) & = z_{1,i}z_{1,j}+z_{2,i}z_{2,j},\quad 1\le i\le j\le n.
\end{align*}
\end{prop}

Set $R=\bC[x_{a,b}]_{(a,b)\in D}$. We have a surjection $R\twoheadrightarrow\cA$ taking $x_{a,b}$ to $\Delta_{a,b}$, denote its kernel by $I$. We refer to $X=\Spec\cA$ as the \textit{type C cluster variety}, this is an irreducible affine variety of dimension $2n-1$. It is cut out by the ideal $I$ in the affine space $\bC^D$.

\begin{remark}
Following one of the common conventions, we define the cluster variety as the affine spectrum of the cluster algebra. This terminology is used in the highly relevant source~\cite{AHHL2021} and by a variety of other authors, e.g.,~\cite{Mul2012,LamSpe2022,CouDuc2020}. The resulting affine scheme is, in general, different from the quasi-affine union of cluster tori, known as the \textit{cluster $\cA$-variety} or the \textit{cluster manifold}.    
\end{remark}

In order to state our main results, we recall the notion of tropicalization in the context of our variety $X$; a comprehensive treatment of these subjects can be found in~\cite[Chapters 2--3]{MacStu2015}. A real weight $w \in \bR^D$ can be viewed as an $\bR$-grading on $R$ that takes the value $w_{a,b}$ on $x_{a,b}$. The \textit{initial form} $\init_w r$ of a polynomial $r\in R$ is its highest nonzero $w$-homogeneous component, i.e., that of maximal grading. For an ideal $I \subset R$, its \textit{initial ideal} $\init_w I$ is the ideal spanned by all initial forms $\init_w r$ with $r\in I$.

For any weight $w\in\bR^D$, the set of weights $w'$ such that $\init_{w'} I=\init_w I$ is a relatively open convex polyhedral cone in $\bR^D$, known as a \textit{Gr\"obner cone} of $I$. Such cones form a complete polyhedral fan in $\bR^D$, the \textit{Gröbner fan} of $I$. 

\begin{defn}
If for $w\in\bR^D$ the ideal $\init_w I$ is monomial-free, then the Gr\"obner cone of $I$ that contains $w$ is a \textit{tropical cone} of $I$. Tropical cones form a subfan $\Trop I$ of the Gröbner fan, known as the \textit{tropicalization} of $I$.    
\end{defn}

It is common to use the notations $\Trop I$ and $\Trop X$ interchangeably, and refer to this fan as the tropicalization of the variety $X$, although it depends not only on $X$ but also on the chosen affine embedding. 

We also use the term \textit{type C tropical cluster variety} to refer to $\Trop X = \Trop I$. The central result of this study is as follows. 

\begin{thm}\label{thm:main}
The type C tropical cluster variety $\Trop X$ coincides with the space of ASPTs $\Omega$.
\end{thm}

Now, consider an arbitrary full rank geometric type cluster algebra of type $\rC_{n-1}$. This is the algebra $\cA(\tilde B)$ given by an $(n-1+m)\times(n-1)$ full rank extended exchange matrix $\tilde B$ of cluster type $\rC_{n-1}$, see~\cite[Section 2]{AHHL2021}. Let $D^*\subset D$ consist of all pairs in $D$ other than $(i,i+1)$ with $i\in[1,n-1]$ and $(1,\ol n)$. Then, $\cA(\tilde B)$ is generated by $n(n-1)$ cluster variables indexed by $D^*$, together with $m$ frozen variables and their inverses. Hence, $\cA(\tilde B)$ is the quotient of $\bC[x_1^{\pm1},\dots,x_m^{\pm1}][x_{a,b}]_{(a,b)\in D^*}$ by an ideal $I(\tilde B)$. Similarly to the above, one may consider initial ideals of $I(\tilde B)$ with respect to weights $w\in\bR^{D^*\cup[1,m]}$, and realize the tropicalization of the cluster variety $X(\tilde B)=\Spec\cA(\tilde B)$ as a fan in $\bR^{D^*\cup[1,m]}$. As shown in~\cite{LamSpe2022}, the group $H(\tilde B)$ of all automorphisms of $\cA(\tilde B)$ that act by scaling the cluster and frozen variables is an $m$-dimensional algebraic group. This implies that the lineality space of $\Trop X(\tilde B)=\Trop I(\tilde B)$ has dimension $m$. 

Let $\mr\cA(\tilde B)$ be the localization of $\cA(\tilde B)$ in all cluster variables, and $\mr X(\tilde B)=\Spec\mr\cA(\tilde B)$ be the very affine part of $X(\tilde B)$. 
The action of $H(\tilde B)$ on $X(\tilde B)$ is free on $\mr X(\tilde B)$. By~\cite[Theorem 4.2]{AHHL2021} the quotient $\mr X(\tilde B)/H(\tilde B)=\Spec\mr\cA(\tilde B)^{H(\tilde B)}$ is independent of $\tilde B$. This quotient is an $(n-1)$-dimensional very affine variety $\cM$ known as the \textit{cluster configuration space}. Its coordinate ring $\bC[\cM]$ is generated by $n(n-1)$ distinguished elements $u_{a,b}$, $(a,b)\in D^*$ and their inverses. This realizes $\Trop\cM$ as a fan in $\bR^{D^*}$. Under the natural isomorphism $\bC[\cM]=\mr\cA(\tilde B)^{H(\tilde B)}$, every $u_{a,b}$ is identified with a certain Laurent monomial in the cluster and frozen variables. This provides a linear surjection from $\Trop X(\tilde B)$ to $\Trop\cM$ whose kernel is the lineality space of the former.

The algebra $\cA$ considered before is almost of the form $\cA(\tilde B)$, one only needs to invert the frozen variables. We have $\cA[x_{a,b}^{-1}]_{(a,b)\in D\bs D^*}=\cA(\tilde B_0)$ for a certain $(n-1)\times(2n-1)$ full rank extended exchange matrix $\tilde B_0$, see~\cite[Example 12.12]{FomZel2003}. In particular, $\Trop X=\Trop X(\tilde B_0)$. Combining this with the above discussion, we have the following

\begin{cor}\label{cor:fullrank}
For any full rank extended exchange matrix $\tilde B$ of type $C_{n-1}$, the tropicalization $\Trop X(\tilde B)$ and the space of ASPTs $\Omega$ are linearly equivalent, modulo lineality space.
\end{cor}

\section{Signed tropicalizations and generalized associahedra}

We next show how signed tropicalizations of the cluster variety and of the cluster configuration space can be interpreted in terms of structures originating in Catalan combinatorics: associahedra and cyclohedra. We start with an overview of positive and signed tropicalizations, in the context of the cluster configuration space $\cM$.

As discussed above, $\bC[\cM]$ is generated by $u_{a,b}^{\pm1}, (a,b)\in D^*$. Hence, $\cM$ is defined by an ideal $J\subset\bC[y_{a,b}^{\pm1}]_{(a,b)\in D^*}=S$ and $\Trop\cM=\Trop J\subset\bR^{D^*}$.
Let $S_{>0}$ denote the subsemiring in $S$ consisting of nonzero Laurent polynomials without negative coefficients. The \textit{positive tropicalization} $\Trop_{>0} \cM=\Trop_{>0} J$ is the subfan of $\Trop \cM$ supported on those $w\in|\Trop \cM|$ for which $\init_w J\cap S_{>0}=\varnothing$.
In other words, $w\in|\Trop_{>0} \cM|$ if and only if every element of $\init_w J$ has both positive and negative coefficients. Let $\cM_{>0} = \cM \cap \bR_{>0}^{D^*}$ be the totally positive part of $\cM$. Clearly, if $\cM_{>0}$ is nonempty, then $\Trop_{>0} \cM$ is also nonempty.

A natural generalization of the positive tropicalization is given by signed tropicalizations. Below, a \textit{sign pattern} is a vector with coordinates in $\{1,-1\}$.
\begin{defn}
For a sign pattern $\nu\in\{1,-1\}^{D^*}$, consider the automorphism $\alpha_\nu$ of $S$ taking $y_{a,b}$ to $\nu_{a,b} y_{a,b}$. The \textit{signed tropicalization} $\Trop_\nu \cM$ is the positive tropicalization $\Trop_{>0} \alpha_\nu(J)$.
\end{defn}

Note that $\Trop\alpha_\nu(J)=\Trop \cM$, hence $\Trop_\nu \cM$ is also a subfan of $\Trop \cM$.

We say that $\nu$ \textit{occurs} in $\cM$ if there is a point $p\in \cM(\bR)=\cM\cap\bR^{D^*}$ such that $\sgn p_{a,b}=\nu_{a,b}$ for all $(a,b)\in D^*$. One sees that $\Trop_\nu \cM$ is nonempty if $\nu$ occurs in $\cM$. We now describe the signed tropicalizations corresponding to all sign patterns that occur in $\cM$. 

To that end, we require the notion of \textit{dihedral orderings}: orderings of $N$ considered up to circular shifts and reversals. There are $(2n-1)!/2$ dihedral orderings in total. With a labeling $\ph$ of the polygon $P$ one naturally associates the dihedral ordering $\la(\ph)$ obtained by reading the labels in $\ph$ clockwise or counter-clockwise, starting from any side. Two labelings provide the same dihedral ordering if and only if they can be identified by the natural action of $D_{4n}$ on $P$.

We distinguish \textit{axially symmetric dihedral orderings} (\textit{ASDOs}): those of the form $\la(\ph)$ with $\ph$ axially symmetric. We also say that a labeling $\ph$ is \textit{centrally symmetric} if the sides $\ph(a)$ and $\ph(\ol a)$ are opposite for every $a\in N$. A \textit{centrally symmetric dihedral ordering} (\textit{CSDO}) is that of the form $\la(\ph)$ with $\ph$ centrally symmetric. There are $2^{n-2}n!$ ASDOs and $2^{n-2}(n-1)!$ CSDOs.

One may also consider \textit{centrally symmetric subdivisions} of $P$: those that with every diagonal contain its image under the central symmetry of $P$. 
\begin{defn}
A \textit{centrally symmetric phylogenetic tree} (or \textit{CSPT}) is a phylogenetic tree of the form $\cT_{\Theta,\ph}$ with $\Theta$ and $\ph$ centrally symmetric.    
\end{defn} 
It is easily seen that every CSPT is also an ASPT but the converse is not true. 

\begin{example}
In the third picture in Figure~\ref{fig:subdivs}, $\Theta$ and $\ph$ are centrally symmetric. Hence, we obtain a CSPT which is simultaneously the ASPT in the first picture. However, the ASPT given by the second picture is not a CSPT.

In Figure~\ref{fig:graphn=3}, vertices and edges corresponding to CSPTs are shown slightly thicker.
\end{example}

Next, we say that a phylogenetic tree $(T,\vv)$ is \textit{compatible} with a dihedral ordering $\la$ if there is a subdivision $\Theta$ and a labeling $\ph$ of $P$ such that $(T,\vv)=\cT_{\Theta,\ph}$ and $\la=\la(\ph)$. An observation easily checked by induction is that an ASPT $(T,\vv)$ is compatible with an ASDO $\la$ if and only if $(T,\vv)$ can be realized as $\cT_{\Theta,\ph}$ with $\Theta$ and $\ph$ \textit{axially symmetric} and $\la(\ph)=\la$. Similarly, a CSPT $(T,\vv)$ is compatible with a CSDO $\la$ if and only if $(T,\vv)$ can be realized as $\cT_{\Theta,\ph}$ with $\Theta$ and $\ph$ centrally symmetric and $\la(\ph)=\la$.

The latter shows that CSPTs compatible with a given CSDO $\la$ are enumerated by all centrally symmetric subdivisions of $P$. Recall that the $(n-1)$-dimensional cyclohedron is a polytope whose face lattice consists of centrally symmetric subdivisions of $P$ ordered by refinement. Recall that $\Trop\cM$ is the lineality-space quotient of the space of ASPTs $\Omega$; denote the image of $C_{T,\vv}$ under this quotient by $C^*_{T,\vv}$. The above shows that cones $C^*_{T,\vv}$ for which $(T,\vv)$ is compatible with $\la$ form a subfan of $\Trop\cM$ which is combinatorially equivalent to the dual fan of a cyclohedron. We denote this subfan by $F_\la$. 

Similarly, ASPTs compatible with a given ASDO $\la$ are enumerated by axially symmetric (with respect to $\delta_0$) subdivisions of $P$. Again, denote by $F_\la$ the subfan of $\Trop\cM$ formed by cones $C^*_{T,\vv}$ for which $(T,\vv)$ is compatible with $\la$. In this case, $F_\la$ is combinatorially equivalent to the dual fan of an $(n-1)$-dimensional associahedron, whose face lattice consists of subdivisions of an $(n+2)$-gon ordered by refinement. To see this, one defines a bijection from the set of axially symmetric subdivisions of $P$ to the set of all subdivisions of an $(n+2)$-gon. Loosely speaking, this bijection contracts the $n-1$ vertices on one side of $\delta_0$ into a single vertex; an example with $n=4$ is shown in Figure~\ref{fig:octtohex}. 

\begin{figure}[h!tbp]
\centering

\begin{tikzpicture}[scale=1/2]
    \foreach \i in {0, 1, 2, 3, 4, 5, 6, 7} {
        \coordinate (a\i) at (45*\i:3);
        \filldraw (a\i) circle (2pt);
        \coordinate (l\i) at ($.55*(45*\i:3) + .55*(45*\i+45:3)$);
    }

    \draw (a0)--(a1)--(a2)--(a3)--(a4)--(a5)--(a6)--(a7)--cycle;

    \draw[dashed] ($1.25*(a1)$)--($1.25*(a5)$);
    \node[right] at ($1.15*(a1)$){$\delta_0$};

    \draw (a4)--(a6);
    \draw (a6)--(a1)--(a4);
    \draw (a4)--(a2);
    \draw (a6)--(a0);

    \draw[->] ($1.7*(a0)$)--($2.7*(a0)$);
\end{tikzpicture}
\hspace{5mm}
\begin{tikzpicture}[scale=1/2]

    \foreach \i in {0, 1, 5, 6, 7} {
        \coordinate (a\i) at (45*\i:3);
        \filldraw (a\i) circle (2pt);
    }

    \draw (a5)--(a6)--(a7)--(a0)--(a1);
    
    \coordinate (v) at ($.5*(a3)$);
    \filldraw (v) circle (2pt);

    \draw (v)--(a6);
    \draw (v)--(a1);
    \draw (v)--(a5);
    \draw (a0)--(a6);
    \draw (a1)--(a6);

     \draw[dashed] ($1.25*(a1)$)--($1.25*(a5)$);

    \node[right] at ($1.15*(a1)$){$\delta_0$};
\end{tikzpicture}

\vspace{-12mm}
\caption{}
\label{fig:octtohex}
\end{figure}

\begin{example}
The 2-dimensional associahedron is a pentagon, hence each $F_\la$ with $\la$ an ASDO is a fan over a 5-cycle. One of the 12 corresponding subgraphs is highlighted in red in Figure~\ref{fig:graphn=3}, it corresponds to the ASDO given by the ordering $(1,2,3,\ol 3,\ol 2,\ol 1)$ of $N$. 
The 2-dimensional cyclohedron is a hexagon. One of the 4 subgraphs corresponding to a subfan $F_\la$ with $\la$ a CSDO is highlighted in blue. 
\end{example}

We can now formulate our description of the signed tropicalizations.

\begin{thm}\label{thm:signedtrops}
There are $2^{n-2}(n+1)(n-1)!$ sign patterns that occur in $\cM$. The respective signed tropicalizations of $\cM$ are the subfans $F_\la$ with $\la$ either an ASDO or a CSDO.
\end{thm}


Now recall that for any full rank cluster variety $X(\tilde B)$, the quotient map to $\cM$ is monomial in the coordinates. Hence, we have a map $f_{\tilde B}$ from sign patterns occurring in $X(\tilde B)$ to those in $\cM$. Moreover, $\Trop_{f_{\tilde B}(\nu)}\cM$ is the lineality-space quotient of $\Trop_\nu X(\tilde B)$, for any $\nu$ that occurs in $X(\tilde B)$. We deduce the following.
\begin{cor}\label{cor:signedtrops}
For any sign pattern $\nu$ occurring in $X(\tilde B)$, the signed tropicalization $\Trop_\nu X(\tilde B)$ is, modulo lineality space and up to combinatorial equivalence, the dual fan of an associahedron or of a cyclohedron.
\end{cor}


\begin{remark}
The positive tropicalization $\Trop_{>0} \cM$ is the subfan $F_{\la}$, where $\la$ is the CSDO given by the ordering $(1,\dots,n,\ol 1,\dots,\ol n)$. A similar statement holds for $\Trop_{>0} X(\tilde B)$, in particular, this fan is dual to a cyclohedron (up to lineality). The latter agrees with a well-known conjecture of Speyer and Williams~\cite{SpeWil2005} settled in~\cite{AHHL2021, JaLoSt2021}, which describes positive tropicalizations of full rank finite type cluster varieties in terms of cluster complexes.
\end{remark}


\acknowledgements{The author is thankful to Lara Bossinger, Shelby Cox, Nathan Ilten and Thomas Lam for valuable discussions.}

\printbibliography

\end{document}